\documentclass[10pt]{article}

\usepackage{graphicx}
\usepackage{epsfig}
\usepackage{amsfonts,amsmath,amsthm,url}

\DeclareMathOperator{\diag}{diag}
\DeclareMathOperator{\tr}{trace}

\theoremstyle{definition}
\newtheorem{definition}{Definition}
\newtheorem{theorem}{Theorem}
\newtheorem{lemma}{Lemma}

\newcommand{\R}{\mathbb{R}}

\newcommand{\cG}{\mathcal{G}}

\newcommand{\cS}{\mathcal{S}}
\newcommand{\vect}{\mbox{vec}}

\title{On semidefinite programming relaxations of the traveling salesman problem}
\author{ Etienne de Klerk\thanks{Department of Econometrics and OR, Tilburg University, The Netherlands. e.deklerk@uvt.nl}
\and Dmitrii V.\ Pasechnik\thanks{School of Physical and Mathematical Sciences,
Nanyang Technological University, Singapore. dima@ntu.edu.sg}
\and Renata Sotirov\thanks{Department of Econometrics and OR, Tilburg University, The Netherlands. r.sotirov@uvt.nl}}

\begin{document}
\maketitle

\begin{abstract}
We consider a new semidefinite programming (SDP) relaxation of the symmetric traveling salesman problem (TSP),
that may be obtained via an SDP relaxation of the more general quadratic assignment problem (QAP).
We show that the new relaxation dominates the one in the paper:

[D.\ Cvetkovi\'c, M.\ Cangalovi\'c and V.\ Kova\v{c}evi\'c-Vuj\v{c}i\'c.
{Semidefinite Programming Methods for the Symmetric Traveling Salesman Problem}.
In  {\em Proceedings of the 7th International IPCO Conference on Integer Programming and Combinatorial Optimization},
{1999},
{126--136},
{Springer-Verlag},
 {London, UK}.]

Unlike the bound of Cvetkovi\'c et al., the new SDP bound is not dominated by the Held-Karp linear programming
bound, or vice versa.

%Finally, we show how to extend the SDP relaxation to a $3/2$-approximation algorithm for metric TSP.
\end{abstract}

\noindent {\bf Keywords:} traveling salesman problem, semidefinite
programming, quadratic assignment problem, association schemes

\vspace{0.3cm} \noindent {\bf AMS classification:} 90C22, 20Cxx,
70-08

\vspace{0.3cm} \noindent {\bf JEL code:} C61

\section{Introduction}

The quadratic assignment problem (QAP) may be stated in the following  form:
\begin{equation}
\label{def:QAP}
\min_{X \in \Pi_n} \tr (AXBX^T)
\end{equation}
where $A$ and $B$ are given symmetric $n\times n$ matrices, and $\Pi_n$ is the set of
$n\times n$ permutation matrices.

It is well-known that the QAP contains the symmetric traveling salesman problem (TSP) as
a special case.
To show this, we denote the complete graph on $n$ vertices
with edge lengths (weights) $D_{ij} = D_{ji} >0$ $(i \neq j)$, by $K_n(D)$,
where $D$ is called the matrix of edge lengths (weights).
The TSP
is to find a Hamiltonian circuit of minimum length in $K_n(D)$.
The $n$ vertices are often called \emph{cities}, and the Hamiltonian circuit of minimum length
the \emph{optimal tour}.

To see that TSP is a special case of QAP,
 let $C_1$ denote the adjacency matrix of ${\cal{C}}_n$ (the standard circuit on $n$ vertices):
\[
C_1:= \left[ \begin{array}{cccccc}
0&1 & 0 &   \cdots & 0 & 1 \\
1&0 & 1 & 0 &  \cdots & 0  \\
 0& 1 & 0 & 1 & \ddots &  \vdots\\
\vdots&\ddots & \ddots & \ddots & \ddots& \\
0& &  & & 0 & 1 \\
1&0 & \cdots & 0 & 1 & 0 \end{array} \right].
\]
Now the TSP problem is obtained from the
QAP problem (\ref{def:QAP}) by setting $A = \frac{1}{2}D$ and $B = C_1$.
To see this, note that
every Hamiltonian circuit in a complete graph has adjacency matrix
$
XC_1X^T$ for some $X \in \Pi_n$. Thus we may concisely state the TSP as
\begin{equation}
TSP_{opt} := \min_{X \in \Pi_n} \tr \left(\frac{1}{2}DXC_1X^T\right).
\label{eq:QAP reformulation of TSP}
\end{equation}

%\subsection*{Complexity and approximation results for QAP and TSP}
The symmetric  TSP is NP-hard in the strong sense \cite{complexity general TSP}, and therefore so is the more general QAP.
In the special case where the distance function of the TSP instance satisfies the triangle inequality (metric TSP), there is a celebrated
$3/2$-approximation algorithm due to Christofides \cite{Christofides heuristic}.
It is a long-standing (since 1975) open problem to improve on the $3/2$ constant, since
the strongest negative result is that a $(1+1/219 )$-approximation algorithm is not possible, unless P=NP
 \cite{Papadimitriou Vempala}.

In the case when the distances are Euclidean in fixed dimension (the so-called planar or geometric TSP),
the problem allows a
polynomial-time approximation scheme \cite{Arora PTAS}.
A recent survey of the TSP is given by Schrijver \cite{Schrijver book}, Chapter 58.

\subsection*{Main results and outline of this paper}
In this paper we will consider semidefinite programming (SDP) relaxations of the TSP.
We will introduce a new SDP relaxation of TSP in Section \ref{sec:new relaxation},
that is motivated by the theory of association schemes.
Subsequently, we will show in Section \ref{sec:relation QAP} that the new SDP relaxation coincides with the SDP relaxation for QAP introduced in \cite{Zhao Karisck Rendl Wolko}
when applied to the QAP reformulation of TSP in (\ref{eq:QAP reformulation of TSP}).
Then we will show in Section \ref{sec:Cvetkovich} that the new SDP relaxation dominates the relaxation due to Cvetkovi\'c et al.\ \cite{Cvetkovic et al 1999}.
The relaxation of Cvetkovi\'c et al. is known to be dominated by the Held-Karp linear programming bound \cite{Dantzig TSP,Held-Karp}, but we show in Section \ref{sec:Held-Karp} that the
new SDP bound is not dominated by the Held-Karp bound (or vice versa).

%Finally, we will show that a family of SDP relaxations dominates the
%minimum length $1$-tree lower bound on $TSP_{opt}$ asymptotically (see Sections \ref{matrix tree theorem} and \ref{sec:Held-Karp bound} for details).

%Finally, we will explore a randomized rounding scheme that uses the optimal solution
%of the new SDP relaxation to construct a Hamiltonian circuit of at most $3/2$ times the optimal value for metric TSP.

\subsection*{Notation}
The space of $p\times q$ real matrices is denoted by $\R^{p\times
q}$,  the space of $k\times k$ symmetric matrices is denoted by
$\cS_k$, and the space of $k\times k$ symmetric positive semidefinite
matrices by $\cS^+_k$.
We will sometimes also use the
notation $X \succeq 0$ instead of
$X \in \cS^+_k$, if the order of the matrix is clear from the context. By $\diag(X)$ we
mean the $n$-vector composed of the diagonal entries of $X \in \cS_n$.

We use $I_n$ to denote the identity
matrix of order $n$. Similarly, $J_n$ and $e_n$ denote the
$n \times n$ all-ones matrix and all ones $n$-vector respectively, and
$0_{n\times n}$ is the zero matrix of order $n$.
We will omit the subscript if the order is clear
from the context.

The {\em Kronecker product} $A \otimes B$ of matrices $A \in
{\R}^{p \times q}$ and $B\in {\R}^{r\times s}$ is
defined as the $pr \times qs$ matrix composed of $pq$ blocks of size
$r\times s$, with block $ij$ given by $A_{ij}B$ $(i = 1,\ldots,p)$,
$(j = 1,\ldots,q)$.

The Hadamard (component-wise) product of matrices $A$ and $B$ of the same size will be denoted by $A \circ B$.

\section{A new SDP relaxation of TSP}
\label{sec:new relaxation}
In this section we show that the optimal value of the following semidefinite program
provides a lower bound on the length $TSP_{opt}$ of an optimal tour:
\begin{equation}
\left.
\begin{array}{rcl}
\min && \frac{1}{2}\tr \left(DX^{(1)}\right) \\
\mbox{subject to}&&\\
&&X^{(k)}  \ge   0, \;\;\;\quad\quad\quad\quad\quad\quad\quad\quad\quad\quad k = 1,\ldots,d \\
&&\sum_{k=1}^d X^{(k)}  =  J-I, \\
&&I + \sum_{k=1}^d \cos\left(\frac{2\pi i k }{n}\right)X^{(k)}   \succeq  0,  \quad\quad i = 1,\ldots, d \\
&&X^{(k)}  \in  \cS^n,  \quad\quad\quad\quad\quad\quad\quad\quad\quad \quad k = 1,\ldots,d,
\end{array}
\right\}
\label{final qap relaxation of TSP}
\end{equation}
where $d = \lfloor \frac{1}{2}n\rfloor$ is the diameter of ${\cal{C}}_n$.

Note that this problem  involves nonnegative matrix variables $X^{(1)},\ldots,X^{(d)}$  of order $n$.
% as opposed to the matrix variable
%of order $n^2$ in  (\ref{final SDP formulation}), i.e.\ the problem size is reduced by a factor $n$ in this sense.
The matrix variables $X^{(k)}$  have an interesting interpretation
in terms of \emph{association schemes}.

\subsection*{Association schemes}
We will give a brief overview of this topic;
 for an introduction to association schemes, see Chapter 12 in \cite{Godsil association schemes}, and in the
context of SDP, \cite{Goemans-Rendl association schemes}.

\begin{definition}[Asssociation scheme]
Assume that a given set of $n\times n$ matrices $B_0,\ldots, B_t$ has the
following properties:
\begin{itemize}
\item[(1)] $B_i$ is a $0-1$ matrix for all $i$ and $B_0 = I$;
\item[(2)] $\sum_i
B_i = J$;
\item[(3)]  $B_i = B_{i^*}^T$
for some $i^*$;
\item[(4)]
$B_iB_j = B_jB_i$ for all $i,j$;
\item[(5)]
$B_iB_j \in \mbox{span}\{B_1,\ldots, B_t\}$.
\end{itemize}
Then we refer to $\{B_1,\ldots, B_t\}$ as an association scheme.
If the $B_i$'s are also symmetric, then we speak of a symmetric association scheme.
\end{definition}
Note that item (4) (commutativity) implies that the matrices $B_1,\ldots, B_t$ share a
common set of eigenvectors, and therefore can be simultaneously diagonalized.
Note also that an association scheme is a basis of a matrix-$*$ algebra (viewed as a vector space).
Moreover, one clearly has
\[
\tr (B_iB_j^T) = 0 \mbox{ if $i \neq j$.}
\]
Since the $B_i$'s share a system of eigenvectors, there is a natural ordering
of their eigenvalues with respect to any fixed ordering of the eigenvectors.
Thus the last equality may be interpreted as:
\begin{equation}
\label{orthogonal eigenvalues}
\sum_k \lambda_k(B_i)\lambda_k(B_j)  = 0 \mbox{ if $i \neq j$,}
\end{equation}
where the $\lambda_k(B_i)$'s  are the eigenvalues of $B_i$ with respect to the fixed ordering.

The association scheme of particular interest to us arises as follows.
Given a connected graph $G = (V,E)$ with diameter $d$, we define
$|V|\times |V|$ matrices $A^{(k)}$ $(k=1,\ldots,d)$ as follows:
\[
A^{(k)}_{ij} = \left\{
\begin{array}{cl}
1 & \mbox{if $\mbox{dist}(i,j) = k$} \\
0 & \mbox{else},
\end{array}
\right.
\quad \quad (i,j \in V),
\]
where dist$(i,j)$ is the length of the shortest path from $i$ to $j$.

Note that $A^{(1)}$ is simply the adjacency matrix of $G$. Moreover, one clearly has
\[
I + \sum_{k=1}^d A^{(k)} = J.
\]
It is well-known  that, for $G={\cal{C}}_n$,  the matrices $A^{(k)}$ $(k=1,\ldots,d\equiv \lfloor n/2 \rfloor)$ together with $A^{(0)} :=I$
form an association scheme, since ${\cal{C}}_n$ is a distance regular graph.

It is shown in the Appendix to this paper, that for $G={\cal{C}}_n$,
the eigenvalues of the matrix $A^{(k)}$ are:
\[
\lambda_m(A^{(k)}) =  2\cos(2\pi mk/n), \quad m = 0,\ldots,n-1, \; k = 1,\ldots,\lfloor (n-1)/2 \rfloor,
\]
and, if $n$ is even,
\[
\lambda_{n/2}(A^{(k)}) = \cos(k\pi) = (-1)^k.
\]
In particular, we have
\begin{equation}
\lambda_m(A^{(k)}) = \lambda_k(A^{(m)}) \quad\quad k,m = 1,\ldots, \lfloor (n-1)/2 \rfloor.
\label{eigenvalue symmetry}
\end{equation}
Also note that
\begin{equation}
\label{eigenvalue symmetry 2}
\lambda_m(A^{(k)}) = \lambda_{n-m}(A^{(k)}), \quad k,m = 1,\ldots, \lfloor (n-1)/2 \rfloor,
\end{equation}
so that each matrix $A^{(k)}$ $(k = 1,\ldots,d)$ has only $1+ \lfloor n/2 \rfloor$ distinct eigenvalues.

\subsection*{Verifying the SDP relaxation (\ref{final qap relaxation of TSP})}
We now show that setting $X^{(k)} = A^{(k)}$ $(k = 1,\ldots,d)$ gives
a feasible solution of (\ref{final qap relaxation of TSP}).
We only need to verify that
\[
I + \sum_{k=1}^d \cos\left(\frac{2\pi i k }{n}\right)A^{(k)}  \succeq  0, \;\; \quad\quad i = 1,\ldots, d.
\]
We will show this for odd $n$, the proof for even $n$ being similar.

Since the $A^{(k)}$'s may be simultaneously diagonalized, the last LMI is the same as:
\[
2 + \sum_{k=1}^d \lambda_k (A^{(i)})\lambda_j (A^{(k)})  \ge   0, \;\; \quad\quad i,j = 1,\ldots, d,
\]
and by using (\ref{eigenvalue symmetry}) this becomes:
\[
2 + \sum_{k=1}^d \lambda_k (A^{(i)})\lambda_k (A^{(j)})  \ge   0, \;\; \quad\quad i,j = 1,\ldots, d.
\]
Since $\lambda_0(A^{(i)}) = 2$ $(i = 1,\ldots, d)$, and using (\ref{orthogonal eigenvalues}), one can easily
verify that the last inequality  holds. Indeed,  one has
\begin{eqnarray*}
&& 2 + \sum_{k=1}^d \lambda_k (A^{(i)})\lambda_k (A^{(j)})\\
&=& 2 + \frac{1}{2} \sum_{k=1}^{n-1} \lambda_k (A^{(i)})\lambda_k (A^{(j)}) \quad \mbox{(by (\ref{eigenvalue symmetry 2}))}\\
&=& 2 - \frac{1}{2}\lambda_0 (A^{(i)})\lambda_0 (A^{(j)}) +  \frac{1}{2} \sum_{k=0}^{n-1} \lambda_k (A^{(i)})\lambda_k (A^{(j)}) \\
&=&
\left\{
\begin{array}{cl}
2 - 2 + 0 = 0 & \mbox{ if $(i\ne j)$, by (\ref{orthogonal eigenvalues})}\\
2-2+\frac{1}{2} \sum_{k=0}^{n-1} \left(\lambda_k (A^{(i)})\right)^2 \ge 0&  \mbox{ if $(i= j)$}.
\end{array}
\right.
\end{eqnarray*}

Thus we have established the following result.

\begin{theorem}
The optimal value of the SDP problem (\ref{final qap relaxation of TSP}) provides a
lower bound on the optimal value $TSP_{opt}$ of the associated TSP instance.
\end{theorem}

\section{Relation of (\ref{final qap relaxation of TSP}) to an SDP relaxation of QAP}
\label{sec:relation QAP}
An SDP relaxation of the QAP problem (\ref{def:QAP})
was introduced in \cite{Zhao Karisck Rendl Wolko}, and
further studied for specially structured instances in \cite{Klerk-Sotirov QAP}.

When applied to the QAP reformulation of TSP in (\ref{eq:QAP reformulation of TSP}), this SDP relaxation takes the form:

\begin{equation}
\label{final SDP formulation}
\left.
\begin{array}{rcl}
\min  & & \frac{1}{2}\tr (C_1\otimes D)Y  \\
\mbox{subject to} && \\
&&\tr((I \otimes (J-I))Y +
((J-I)\otimes I) Y) = 0 \\
&&\tr (Y) - 2e^Ty = -n \\
% Y_{00} &=& 1 \\
&&\left(
\begin{array}{cc}
1 & y^T \\
y & Y
\end{array}
\right) \succeq 0, \;\;
Y \geq  0. \\
\end{array}
\right\}
\end{equation}

It is easy to verify that this is indeed a relaxation of problem (\ref{eq:QAP reformulation of TSP}),
by noting that setting
$Y = \vect(X)\vect(X)^T$ and  $y = \diag(Y)$ gives  a feasible solution if $X \in \Pi_n$.

In this section we will show that the optimal value of the SDP problem (\ref{final SDP formulation}) actually
 equals the optimal value of the new SDP relaxation (\ref{final qap relaxation of TSP}).
The proof is via the technique of \emph{symmetry reduction}.

\subsection*{Symmetry reduction of the SDP problem (\ref{final SDP formulation})}
Consider the following form  of  a general semidefinite programming problem:
\begin{equation} \label{sdp}
p^* := \min\limits_{X\succeq 0, X \ge 0} \left \{ ~\tr(A_0 X) \; : \;
\tr(A_kX)=b_k,\quad k=1,\ldots,m \right \},
\end{equation}
where the $A_i$ $(i=0,\ldots,m)$ are given symmetric matrices.

If we view (\ref{final SDP formulation}) as an SDP problem in the  form (\ref{sdp}),
 the data matrices of problem
(\ref{final SDP formulation})
are:
\begin{equation}
\small
\label{SDP TSP data}
\left(
\begin{array}{cc}
0 & 0^T \\
0 & \frac{1}{2}C_1 \otimes D
\end{array}
\right),
\left(
\begin{array}{cc}
0 & 0^T \\
0 & I \otimes (J-I) +
(J-I)\otimes I,
\end{array}
\right),
\left(
\begin{array}{cc}
0 & -e^T \\
-e & 2I
\end{array}
\right),
\left(
\begin{array}{cc}
1 & 0^T \\
0 & 0
\end{array}
\right).
\end{equation}

\begin{definition}
We define the \emph{automorphism group} of a  matrix $Z \in \mathbb{R}^{k\times k}$ as
\[
\mbox{aut}(Z) = \{P \in \Pi_k \; : \; PZP^T = Z\}.
\]
\end{definition}

Symmetry reduction of problem (\ref{sdp}) is possible under the  assumption that
the  multiplicative matrix group
\[
 {\cal G} := \bigcap_{i=0}^m \mbox{aut}(A_i)
\]
is non-trivial. We call ${\cal G}$ the symmetry group of the SDP problem (\ref{sdp}).

For the matrices (\ref{SDP TSP data}), the group
${\cal G}$ is given by the matrices
\begin{equation}
\label{tsp group}
\cG := \left\{
\left(
\begin{array}{cc}
1 & 0^T \\
0 & P \otimes I
\end{array}
\right) \; : \; P  \in \mathcal{D}_n \right\},
\end{equation}
where $\mathcal{D}_n$ is the (permutation matrix representation of) the dihedral group of order $n$, i.e.\ the automorphism group of  ${\cal{C}}_n$.

The basic idea of symmetry reduction is given by the following result.
\begin{theorem}[see e.g.\ \cite{GatPa:04}]
If $X$ is a feasible (resp.\ optimal) solution of the SDP problem (\ref{sdp}) with symmetry group ${\cal G}$,
then
\[
\bar X := \frac{1}{|{\cal G}|}\sum_{P\in {\cal G}}
P^TXP
\]
is also a feasible (resp.\ optimal) solution of (\ref{sdp}).
\end{theorem}

Thus there exist  optimal solutions in the set
\[
{\cal{A}}_{\cal{G}} := \left\{ {\frac{1}{|{\cal G}|}\sum_{P\in {\cal G}} P^TXP} \; : \; X \in \mathbb{R}^{  n\times  n} \right\}.
\]
This set is called the centralizer ring (or commutant) of ${\cal G}$ and it is a matrix $*$-algebra.
For the group defined in
(\ref{tsp group}), it is straightforward to verify that the centralizer ring is given by:
\begin{equation}
{\cal{A}}_{\cal{G}}:=\left. \left\{
\left(
\begin{array}{cc}
\alpha & x^T \\
y & C\otimes Z
\end{array}
\right) \; \right| \;
\alpha \in \mathbb{R}, \; C=C^T \mbox{ circulant}, \; Z \in \mathbb{R}^{n\times n}, \; x,y \in \mathbb{R}^{n^2}
\right\}
\label{block form of feasible solutions}
\end{equation}
where $x^T = [x_1e^T \ldots x_ne^T]$ and $y^T = [y_1e^T \ldots y_ne^T]$ for some
scalars $x_i$ and $y_i$ $(i=1,\ldots,n)$, where $e \in \mathbb{R}^n$ is the all-ones vector, as before.

Thus we may restrict the feasible set of problem (\ref{final SDP formulation}) to feasible solutions of the form
(\ref{block form of feasible solutions}).

If we divide $y$ and $Y$ in (\ref{final SDP formulation}) into blocks:
\[
y = \left(\left(y^{(1)}\right)^T \cdots \left(y^{(n)}\right)^T\right)^T,
\]
and
\[
Y =
\left(
\begin{array}{ccc}
Y^{(11)} & \cdots & Y^{(1n)} \\
\vdots &\ddots & \vdots \\
Y^{(n1)} &\cdots& Y^{(nn)}
\end{array}\right),
\]
where $y^{(i)} \in \mathbb{R}^{n}$ and $Y^{(ij)} = Y^{(ji)T} \in\mathbb{R}^{n\times n}$,
then feasible solutions of (\ref{final SDP formulation}) satisfy
\begin{equation}
\label{block form}
 \left(
\begin{array}{cccc}
1 & \left(y^{(1)}\right)^T & \cdots & \left(y^{(n)}\right)^T \\
y^{(1)}& Y^{(11)} & \cdots & Y^{(1n)} \\
\vdots & \vdots &\ddots & \vdots \\
y^{(n)} & Y^{(n1)} &\cdots& Y^{(nn)}
\end{array}\right) \succeq 0.
\end{equation}

Feasible solutions have the following additional structure (see \cite{Zhao Karisck Rendl Wolko} and Theorem 3.1 in \cite{Klerk-Sotirov QAP}):
\begin{itemize}
\item
$Y^{(ii)}$ $(i=1,\ldots,n)$ is a diagonal matrix;
\item
 $Y^{(ij)}$ $(i \neq j)$
is a matrix with zero diagonal;
\item
  $\tr(JY^{(ij)}) = 1$ $(i,j = 1,\ldots,n)$;
  \item
  $\sum_{i=1}^n Y^{(ij)} = e\left(y^{(j)}\right)^T \quad (j=1,\ldots,n)$;
\item
$\diag(Y) = y$.
\end{itemize}

Since $\diag(Y) = y$ for feasible solutions, we have $y^{(i)} = \diag(Y^{(ii)})$ $(i=1,\ldots,n)$.
Moreover, since we may also assume the structure (\ref{block form of feasible solutions}),
we have that
\[
y^{(i)} = y_ie \quad (i=1,\ldots,n),
\]
for some scalar values $y_i$. This implies that the diagonal elements of $Y^{(ii)}$ all equal $y_i$.
Since the diagonal elements of $Y^{(ii)}$ sum to $1$, we have $y_i = 1/n$ and $\diag(Y^{(ii)}) = (1/n)e$.
Thus the condition:
\[
\left(
\begin{array}{cc}
1 & y^T \\
y & Y
\end{array}
\right) \succeq 0
\]
reduces to
\[
Y - \frac{1}{n^2}J \succeq 0
\]
 by the Shur complement theorem.
This is equivalent to
\[
(I \otimes Q^*)Y(I\otimes Q) - \frac{1}{n^2}(I \otimes Q^*)J(I \otimes Q) \succeq 0
\]
where $Q$ is the discrete Fourier transform matrix defined in (\ref{def:discrete Fourier transform matrix}) in the Appendix.

Using the properties of the Kronecker product and of $Q$ we get
\[
\left(
\begin{array}{ccc}
Q^*Y^{(11)}Q & \cdots & Q^*Y^{(1n)}Q \\
\vdots &\ddots & \vdots \\
Q^*Y^{(n1)}Q &\cdots& Q^*Y^{(nn)}Q
\end{array}\right)
-
J\otimes
\left(
\begin{array}{ccc}
\frac{1}{n} & \cdots & 0 \\
\vdots &\ddots & \vdots \\
0 &\cdots& 0
\end{array}\right)
\succeq 0.
\]

Recall that $Y^{(ii)} = \frac{1}{n}I$ and that we may assume
$Y^{(ij)}$ $(i \neq j)$ to be symmetric circulant, say
\[
Y^{(ij)} = \sum_{k=1}^d x^{(ij)}_kC_k,  \quad (i\neq j),
\]
where $C_k$ $(k = 1,\ldots,d)$
forms a basis of the symmetric circulant matrices with zero diagonals (see the Appendix for the precise definition).
Note that the nonnegativity of $Y^{(ij)}$ is equivalent to $x^{(ij)}_k \ge 0$ $(k = 1,\ldots,d)$.
Since $\tr(JY^{(ij)}) = 1$ one has
\[
\sum_{k=1}^d x^{(ij)}_k = \frac{1}{2n} \quad (i\neq j).
\]
Since $\sum_{i=1}^n Y^{(ij)} = e\left(y^{(j)}\right)^T = \frac{1}{n}J$, one  also has
\[
\sum_{k=1}^d \sum_{i=1}^n x^{(ij)}_kC_k = \frac{1}{n}J.
\]
By the definition of the $C_k$'s, this implies that
\begin{equation}
\label{eq:block sums}
\sum_{i=1}^n x^{(ij)}_k =
\left\{
\begin{array}{cc}
\frac{1}{n} & \mbox{if $1 \le k \le \lfloor (n-1)/2 \rfloor$} \\
\frac{1}{2n} & \mbox{if $k = n/2$ ($n$ even)}.
\end{array}
\right.
\end{equation}

Moreover,
\[
Q^*Y^{(ij)}Q = \sum_{k=1}^d x^{(ij)}_kD_k,  \quad (i\neq j),
\]
where $D_k$ is the diagonal matrix with the eigenvalues (\ref{eigs of circulant basis}) of $C_k$ on its diagonal.

Thus the LMI becomes
\begin{equation}
\left(
\begin{array}{ccc}
\frac{1}{n}I & \cdots & \sum_{k=1}^d x^{(1n)}_kD_k \\
\vdots &\ddots & \vdots \\
\sum_{k=1}^d x^{(1n)}_kD_k &\cdots& \frac{1}{n}I
\end{array}\right)
-
J\otimes
\left(
\begin{array}{ccc}
\frac{1}{n} & \cdots & 0 \\
\vdots &\ddots & \vdots \\
0 &\cdots& 0
\end{array}\right)
\succeq 0.
\label{almost final LMI}
\end{equation}
The left hand side of this LMI is a block matrix with each block being a diagonal matrix.
Thus this matrix has a chordal sparsity structure ($n$ disjoint cliques of size $n$).
We may now use the following lemma to obtain the system of LMI's
(\ref{final qap relaxation of TSP}).

\begin{lemma}[cf.\ \cite{chordal}]
Assume a $nt \times nt$ matrix has the block structure
\[
M:= \left(
\begin{array}{ccc}
D^{(11)} & \cdots & D^{(1n)} \\
\vdots &\ddots & \vdots \\
D^{(n1)} &\cdots& D^{(nn)}
\end{array}\right),
\]
where $D^{(ij)} \in \cS_t$ are diagonal $(i,j = 1,\ldots,n)$.
Then $M\succeq 0$ if and only if:
\[
\left(
\begin{array}{ccc}
D^{(11)}_{ii} & \cdots & D^{(1n)}_{ii} \\
\vdots &\ddots & \vdots \\
D^{(n1)}_{ii} &\cdots& D^{(nn)}_{ii}
\end{array}\right) \succeq 0 \quad i = 1,\ldots,t.
\]
\end{lemma}

Applying the lemma to the LMI (\ref{almost final LMI}), and
setting
\begin{equation}
\label{eq:block sums2}
X^{(k)}_{ij} = 2n x_k^{(ij)}, \quad \quad k = 1,\ldots,\lfloor n/2 \rfloor,
\end{equation}
 yields the system of LMI's in (\ref{final qap relaxation of TSP}).

Thus we have established the following result.

\begin{theorem}
The optimal values of the semidefinite programs (\ref{final qap relaxation of TSP}) and (\ref{final SDP formulation})  are equal.
\end{theorem}

\section{Relation of (\ref{final qap relaxation of TSP}) to an SDP relaxation of Cvetkovi\'c et al.}
\label{sec:Cvetkovich}
We will now show that the new SDP relaxation (\ref{final qap relaxation of TSP})  dominates an SDP relaxation
(\ref{Cvetkovich relaxation of TSP}) due to Cvetkovi\'c et al.\ \cite{Cvetkovic et al 1999}.
This latter relaxation is based on the fact that the
spectrum of  the Hamiltonian circuit ${\cal C}_n$ is known.
In particular, the smallest eigenvalue of its Laplacian is zero and corresponds to the all ones eigenvector, while the
second smallest eigenvalue equals $2-2\cos\left(\frac{2\pi}{n} \right)$.

The relaxation takes the form:
\[
TSP_{opt} \ge \min \frac{1}{2} \tr (DX)
\]
subject to
\begin{equation}
\left.
\begin{array}{rcl}
Xe & = & 2e, \\
\mbox{diag}(X) & = & 0, \\
0 & \le & X \le J, \\
2I-X + \left(2-2\cos\left(\frac{2\pi}{n}\right)\right)(J-I) & \succeq & 0.
\end{array}
\right\}
\label{Cvetkovich relaxation of TSP}
\end{equation}
Note that the matrix variable  $X$ corresponds to the adjacency
matrix of the minimal length Hamiltonian circuit.

\begin{theorem}
The SDP relaxation (\ref{final qap relaxation of TSP}) dominates the relaxation (\ref{Cvetkovich relaxation of TSP}).
\end{theorem}
\proof
Assume that given $X^{(k)}  \in  \cS^n$ $(k = 1,\ldots,d)$ satisfy (\ref{final qap relaxation of TSP}).
Then, $$\mbox{diag}(X^{(1)})  =  0,$$
while
(\ref{eq:block sums}) and (\ref{eq:block sums2}) imply
\[
X^{(k)}e = 2e \quad (k = 1,\ldots, \lfloor (n-1)/2 \rfloor),
\]
and
$X^{(n/2)}e = e$ if $n$ is even.
In particular, one has
$X^{(1)}e  =  2e$.
It remains to show that
\[
2I-X^{(1)} + \left(2-2\cos\left(\frac{2\pi}{n}\right)\right)(J-I)  \succeq  0,
\]
which is the same as showing that
\begin{equation}
\label{lmi to be derived}
2I-X^{(1)} + \left(2-2\cos\left(\frac{2\pi}{n}\right)\right)\sum_{k=1}^d X^{(k)}  \succeq  0,
\end{equation}
since
\[
\sum_{k=1}^d X^{(k)} = J-I.
\]
We will show that the LMI (\ref{lmi to be derived}) may be obtained as a nonnegative aggregation
of the LMI's
\[
I+  \sum_{k=1}^d X^{(k)}  \succeq 0
\]
and
\[
I + \sum_{k=1}^d \cos\left(\frac{2\pi i k }{n}\right)X^{(k)}   \succeq  0 \quad (i = 1,\ldots,d).
\]
The matrix of coefficients of these LMI's
is a $(d+1) \times (d+1)$ matrix, say $A$, with entries:
\[
A_{ij} = \cos\left(\frac{2\pi i j}{n}\right) \quad (i,j = 0,\ldots,d).
\]
Since we may rewrite (\ref{lmi to be derived}) as
\[
2I +  \left(1-2\cos\left(\frac{2\pi}{n}\right)\right)X^{(1)} + \left(2-2\cos\left(\frac{2\pi}{n}\right)\right)\sum_{k=2}^d X^{(k)}  \succeq  0,
\]
we need to show that the linear system $Ax = b$ has a nonnegative solution,
where
\[
b := \left[2, \left(1-2\cos\left(\frac{2\pi}{n}\right)\right), \left(2-2\cos\left(\frac{2\pi}{n}\right)\right), \ldots, \left(2-2\cos\left(\frac{2\pi}{n}\right)\right)\right]^T.
\]
One may verify that, for $n$ odd, the system $Ax=b$ has a (unique) solution given by
\[
x_i =
\frac{4}{n}\left\{ \begin{array}{ll}
d\left(1-\cos\left(\frac{2\pi}{n}\right)\right) & \mbox{if $i=0$} \\
\cos\left(\frac{2\pi}{n}\right)-\cos\left(\frac{2\pi i}{n}\right) & \mbox{for $i=1,\ldots,d$.}
\end{array}
\right.
\]
Note that $x$ is nonnegative, as it should be. If $n$ is even, the solution is
\[
x_i =
\frac{4}{n}\left\{ \begin{array}{ll}
\frac{(n-1)}{2}\left(1-\cos\left(\frac{2\pi}{n}\right)\right) & \mbox{if $i=0$} \\
\cos\left(\frac{2\pi}{n}\right)-\cos\left(\frac{2\pi i}{n}\right) & \mbox{for $i=1,\ldots,d-1$} \\
\frac{1}{2}\cos\left(\frac{2\pi}{n}\right)-\frac{1}{2}\cos\left(\frac{2\pi i}{n}\right) & \mbox{for $i=d$.} \\
\end{array}
\right.
\]
\qed

In the section with numerical examples, we will present instances where the new SDP
relaxation (\ref{final qap relaxation of TSP}) is
 strictly better than (\ref{Cvetkovich relaxation of TSP}).

\section{Relation to the Held-Karp bound}
\label{sec:Held-Karp}
One of the best-known linear programming (LP) relaxations
of TSP is the  LP with sub-tour elimination constraints:
\[
TSP_{opt} \ge \min \frac{1}{2} \tr (DX)
\]
subject to
\begin{equation}
\left.
\begin{array}{rcl}
Xe & = & 2, \\
\mbox{diag}(X) & = & 0, \\
0 & \le & X \le J, \\
\sum_{i \in {\cal{I}}, \; j \notin {\cal{I}}} X_{ij} & \ge & 2 \quad \forall  \;   \emptyset \neq   {\cal{I}} \subset \{1,\ldots,n\}.
\end{array}
\right\}
\label{subtour elimination relaxation of TSP}
\end{equation}
This LP relaxation dates back to 1954 and is due to Dantzig, Fulkerson and Johnson \cite{Dantzig TSP}. Its optimal value coincides with the LP bound of
Held and Karp \cite{Held-Karp} (see e.g.\ Theorem 21.34 in \cite{Korte-Vygen}), and the optimal value of the LP is commonly known as the \emph{Held-Karp bound}.

The last constraints are called \emph{sub-tour elimination inequalities} and model the fact that ${\cal{C}}_n$ is 2-connected.
Although there are exponentially many sub-tour elimination inequalities, it is well-known that
the LP (\ref{subtour elimination relaxation of TSP}) may be solved in polynomial time using the ellipsoid method; see e.g.\
Schrijver \cite{Schrijver book}, \S 58.5. %Indeed, the ellipsoid method only requires a polynomial time separation oracle
%for infeasible $X$, and a suitable separating hyperplane is obtained by
%solving the minimum cut problem on the graph with capacities matrix $X$.

It was shown by Goemans and Rendl \cite{Goemans-Rendl} that this LP relaxation dominates the SDP relaxation (\ref{Cvetkovich relaxation of TSP})
 by Cvetkovi\'c et al. \cite{Cvetkovic et al 1999}.
The next theorem shows that
the LP relaxation (\ref{subtour elimination relaxation of TSP}) does not dominate the new SDP relaxation
(\ref{final qap relaxation of TSP}), or vice versa.

\begin{theorem}
\label{th:no dominating relaxation}
The LP sub-tour elimination relaxation (\ref{subtour elimination relaxation of TSP}) does not dominate the new SDP relaxation
(\ref{final qap relaxation of TSP}), or vice versa.
\end{theorem}
\proof

\begin{figure}[h!]
\begin{center}
\includegraphics{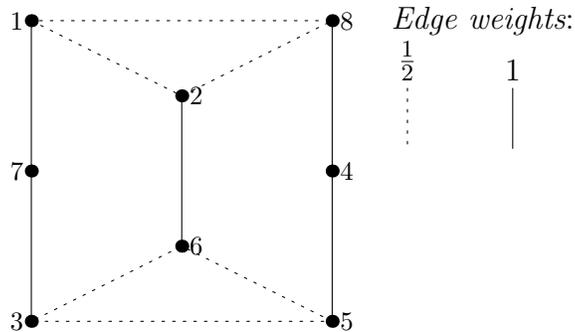}
\caption{The weighted graph used in the proof of Theorem \ref{th:no dominating relaxation}. \label{fig:barX}}
\end{center}
\end{figure}

Define the $8 \times 8$ symmetric matrix $\bar X$ as the weighted adjacency matrix of the graph shown  in Figure \ref{fig:barX}.

The matrix $\bar X$ satisfies the sub-tour elimination inequalities, since the minimum cut in the graph in
Figure \ref{fig:barX} has weight $2$.

On the other hand, there does not exist a feasible solution of (\ref{final qap relaxation of TSP})
that satisfies $X^{(1)} = \bar X$, as may be shown using SDP duality theory.

Conversely, in Section \ref{sec:numerical examples} we will provide  examples where the optimal value
of (\ref{subtour elimination relaxation of TSP}) is strictly greater than  the optimal value
of (\ref{final qap relaxation of TSP}) (see e.g.\ the instances gr17, gr24 and bays24 there). \qed

\section{An LMI cut via the number of spanning trees}
\label{sec:matrix tree cut}
In addition to the sub-tour elimination inequalities, there are several families of
 linear inequalities known for the TSP polytope;
for a review, see Naddef \cite{Naddef} and Schrijver \cite{Schrijver book}, Chapter 58.

Of particular interest to us is a valid nonlinear inequality that models the fact that
${\cal{C}}_n$ has $n$ distinct spanning trees.
 To introduce the inequality we require
a general form of the {\em matrix tree theorem}; see {\em e.g.}\ Theorem VI.29 in \cite{Tutte} for a proof.

\begin{theorem}[Matrix tree theorem]
\label{matrix tree theorem}
Let a simple graph $G = (V,E)$ be given and associate with each edge $e \in E$
a real variable $x_e$.
Define the (generalized) Laplacian of $G$ with respect to $x$ as
the $|V|\times |V|$ matrix with entries
\[
L(G)(x)_{ij} :=
\left\{
\begin{array}{ll}
\sum_{e \; : \; e \cap i \neq \emptyset}x_e & \mbox{if ${i=j}$}, \\
-x_e & \mbox{if $\{i,j\}=e$}, \\
0 & \mbox{else}.
\end{array}
\right.
\]
Now all principal minors of $L(G)(x)$ of order $|V|-1$ equal:
\begin{equation}
\label{number of spanning trees}
\sum_{T} \prod_{e \in T} x_e,
\end{equation}
where the  sum is over all distinct spanning trees $T$ of $G$.
\end{theorem}

In particular, if $L(G)(x)$ is the usual Laplacian of a given graph, then $x_e = 1$ for all edges $e$ of the graph,
and expression (\ref{number of spanning trees}) evaluates to the number of spanning trees in the graph.

Thus if $X$ corresponds to the approximation
of the adjacency matrix of a minimum tour, then
one may require that:
\begin{equation}
\det\left(2I-X\right)_{2:n,2:n} \ge n,
\label{minor inequality}
\end{equation}
where $X_{2:n,2:n}$ denotes the principle submatrix of $X$ obtained by deleting the first row and column.

The inequality (\ref{minor inequality}) may be added to the above SDP relaxations (\ref{Cvetkovich relaxation of TSP})
and (\ref{final qap relaxation of TSP}) (with $X = X^{(1)}$),
since the set
\[
\{ Z \succeq 0 \; : \; \det Z \ge n\}
\]
is LMI representable; see e.g.\ Nemirovski \cite{Nemirovskii lecture notes}, \S3.2.

We know from numerical examples that (\ref{minor inequality}) is not implied
by the relaxation of Cvetkovi\'c et al.\ (\ref{Cvetkovich relaxation of TSP}), but do not know
any examples where it is violated by a feasible $X^{(1)}$
of the new relaxation (\ref{final qap relaxation of TSP}). Nevertheless,
we have been unable to show that (\ref{minor inequality}) (with $X = X^{(1)}$) is implied by
(\ref{final qap relaxation of TSP}).

\section{Numerical examples}
\label{sec:numerical examples}
In Table \ref{tab:bounds} we give the
lower bounds on some small TSPLIB\footnote{\url{http://www.iwr.uni-heidelberg.de/groups/comopt/software/TSPLIB95/}}
 instances for the two SDP relaxations (\ref{final qap relaxation of TSP}) and  (\ref{Cvetkovich relaxation of TSP}), as
 well as the LP relaxation with all sub-tour elimination constraints (\ref{subtour elimination relaxation of TSP}) (the Held-Karp bound).
 These instances have integer data, and the optimal values of the relaxations were rounded up to obtain the bounds in the table.

The  SDP problems were solved by the interior point software
CSDP \cite{CSDP} using the Yalmip interface \cite{YALMIP}
and Matlab 6.5, running on a PC with two 2.1 GHz dual-core processors and 2GB of  memory.

\begin{table}[h!]
{\small
\begin{center}
\begin{tabular}{|c|c|c|c|c|} \hline
Problem & SDP bound (\ref{Cvetkovich relaxation of TSP})  & SDP bound (\ref{final qap relaxation of TSP}) (time) &  LP bound (\ref{subtour elimination relaxation of TSP})& $TSP_{opt}$ \\ \hline \hline
gr17    &  1810     & 2007 (39s) &    2085 & 2085 \\ \hline
gr21    &  2707     & 2707  (139s) &    2707 & 2707 \\ \hline
gr24    &  1230     & 1271  (1046s)&   1272 & 1272  \\ \hline
bays29  &  1948     & 2000 (2863s) &   2014 & 2020 \\ \hline
\end{tabular}
\caption{\label{tab:bounds} Lower bounds on some small TSPLIB instances from various convex relaxations.}
\end{center}
}
\end{table}

Note that the relaxation (\ref{final qap relaxation of TSP})
can indeed be strictly better than (\ref{Cvetkovich relaxation of TSP}), as is clear from the gr17, bays24 and bays29 instances.
Also, since the LP relaxation (\ref{subtour elimination relaxation of TSP}) gives better bounds than (\ref{final qap relaxation of TSP}) for all four
instances, it is worth recalling that this will not happen in general, by Theorem
\ref{th:no dominating relaxation}.

The LMI cut from (\ref{minor inequality}) was already satisfied by the optimal solutions of (\ref{Cvetkovich relaxation of TSP}) and
(\ref{final qap relaxation of TSP}) for the four instances.

A second set of test problems was generated by considering
all facet defining inequalities for the TSP polytope on $8$ nodes; see \cite{TSP on 8 nodes} for a description
of these inequalities, as well as  the SMAPO project web site\footnote{\url{http://www.iwr.uni-heidelberg.de/groups/comopt/software/SMAPO/tsp/}}.

The facet defining inequalities are of the form $\frac{1}{2}\tr (DX) \ge RHS$ where
 $D \in {\cal{S}}_n$ has nonnegative integer entries  and $RHS$ is an integer. From each inequality, we form
a symmetric TSP instance with distance matrix $D$. Thus the optimal value of the TSP instance is the
value $RHS$.
In Table \ref{tab:nis8} we give the optimal values of the LP relaxation (\ref{subtour elimination relaxation of TSP}) ({\em i.e.}\ the Held-Karp bound),
the SDP relaxation of Cvetkovi\'c et al.\ (\ref{Cvetkovich relaxation of TSP}),
 and the new SDP relaxation (\ref{final qap relaxation of TSP})
for these instances, as well as the right-hand-side $RHS$ of each  inequality $\frac{1}{2}\tr (DX) \ge RHS$.
For $n=8$, there are $24$ classes of facet-defining inequalities. The members of each class are equal modulo a permutation of the nodes,
and we need therefore only consider one representative per class. The first
 three classes of inequalities are sub-tour elimination
inequalities.

\begin{table}[h!]
\begin{center}
{\footnotesize
\begin{tabular}{|l|l|l|l|l|} \hline
Inequality & SDP bound (\ref{Cvetkovich relaxation of TSP}) & SDP bound (\ref{final qap relaxation of TSP}) &  Held-Karp bound (\ref{subtour elimination relaxation of TSP}) & RHS  \\ \hline \hline
1  &  2  &  2&  2& 2 \\ \hline
2  &  1.098  & 1.628&  2& 2\\ \hline
3  &  1.172  & 1.172& 2& 2 \\ \hline
4  &  8.507  & 8.671 & 9& 10  \\ \hline
5  &  9  & 9 & 9& 10  \\ \hline
6  &  8.566  & 8.926 & 9& 10  \\ \hline
7  &  8.586  & 8.586 & 9& 10  \\ \hline
8  &  8.570  & 8.926&  9& 10  \\ \hline
9  &  9  & 9 & 9& 10  \\ \hline
10 &  8.411  & 8.902&  9& 10  \\ \hline
11 &  8.422  & 8.899&  9& 10  \\ \hline
12 &  0  & 0& 0 & 0  \\ \hline
13 &  10.586 & 10.667 & 11 & 12  \\ \hline
14 &  12 & 12 & 12 & 13 \\ \hline
15 &  12.408 & 12.444 & 12$\frac{2}{3}$ & 14 \\ \hline
16 &  14 & 14.078 & 14 & 16 \\ \hline
17 &  16 & 16 & 16 & 18 \\ \hline
18 &  16 & 16 & 16  & 18 \\ \hline
19 &  16 & 16 & 16 & 18 \\ \hline
20 &  15.185 & 15.926 & 16 & 18 \\ \hline
21 &  18 & 18.025 & 18 & 20 \\ \hline
22 &  20 & 20 & 20 & 22 \\ \hline
23 &  23 & 23.033 & 23 & 26 \\ \hline
24 &  34.586 & 34.739 & 35 & 38  \\ \hline
\end{tabular}}
\caption{Results for instances on $n=8$ cities, constructed from the facet defining inequalities. \label{tab:nis8}}
\end{center}
\end{table}

The numbering of the instances in Table \ref{tab:nis8} coincides with the numbering of the classes of facet defining inequalities on the
SMAPO project web site.

The new SDP bound (\ref{final qap relaxation of TSP}) is only stronger than the Held-Karp bound (\ref{subtour elimination relaxation of TSP})
for the instances 16, 21 and 23 in Table \ref{tab:nis8}, and for the instances 1, 5, 9, 12, 14, 17, 18, 19 and 22 the two bounds coincide.
For the remaining 18 instances the Held-Karp bound is better than the SDP bound (\ref{final qap relaxation of TSP}).
 However, if the bounds are rounded up, the SDP bound (\ref{final qap relaxation of TSP}) is still better for the instances 16, 21 and 23, whereas the
 two (rounded) bounds are
equal for all the other instances.
Adding the LMI cut from (\ref{minor inequality}) did not change the optimal values  of
the SDP relaxations (\ref{Cvetkovich relaxation of TSP}) or
(\ref{final qap relaxation of TSP}) for any of the instances.

For $n = 9$, there are 192 classes of facet defining inequalities of the TSP polytope \cite{TSP on 9 nodes}.
Here the
SDP bound (\ref{final qap relaxation of TSP}) is better than the Held-Karp bound for 23 out of the 192 associated TSP instances.
Similar to the $n=8$ case, when rounding up, the rounded SDP bound remains better in all 23 cases and coincides with the
rounded Held-Karp bound in all the remaining cases.

\section{Concluding remarks}

Wolsey \cite{Wolsey TSP} showed that the optimal value of the LP relaxation (\ref{subtour elimination relaxation of TSP})
is at least $2/3$ the length of an optimal tour for metric TSP (see also \cite{Shmoys-Williamson1990}).
 An interesting question is whether a similar result may be proved
for the new SDP relaxation (\ref{final qap relaxation of TSP}).

\vspace{0.3cm}

Finally, the computational perspectives of the SDP relaxation (\ref{final qap relaxation of TSP})
are somewhat limited due to its size.
However, since it provides a new polynomial-time convex approximation
of TSP with a rich mathematical structure, it is our hope that
it may lead to a renewed interest in improving approximation results for metric TSP.

\subsection*{Acknowledgements}
Etienne de Klerk would like to thank Drago\u{s} Cvetkovi\'c and Vera Kova\v{c}evi\'c-Vuj\v{c}i\'c for past discussions
on the SDP relaxation (\ref{Cvetkovich relaxation of TSP}). The authors would also like to thank
an anonymous referee for suggestions that led to a significant improvement of this paper.

\section*{Appendix: Circulant matrices}
Our discussion of circulant matrices is condensed from
the review paper by Gray \cite{circulant matrices}.

A circulant matrix has the form
\begin{equation}
\label{def:circulant matrix}
C= \left[ \begin{array}{cccccc}
c_0& c_1 & c_2 &   \cdots &  & c_{n-1} \\
c_{n-1} &c_0 & c_1 &  &   &   \\
 & c_{n-1} & c_0 & c_1 &  &  \vdots\\
\vdots&\ddots & \ddots & \ddots & \ddots& \\
& &  & &  & c_1 \\
c_1& & \cdots &  & c_{n-1} & c_0 \end{array} \right].
\end{equation}
Thus the entries satisfy the
relation
\begin{equation}
\label{cyclic symmetry of circulant}
C_{ij} = c_{(j-i) \mod n}.
\end{equation}

The matrix $C$ has eigenvalues
\[
\lambda_m(C) = c_0 + \sum_{k=1}^{n-1} c_ke^{-2\pi \sqrt{-1} mk/n}, \quad m = 0,\ldots,n-1.
\]
If $C$ is symmetric with $n$ odd, this reduces to
\begin{equation}
\label{eigs symmetric circulant}
\lambda_m(C) = c_0 + \sum_{k=1}^{(n-1)/2} 2c_k\cos(2\pi mk/n), \quad m = 0,\ldots,n-1,
\end{equation}
and when $n$ is even we have
\begin{equation}
\label{eigs symmetric circulant even case}
\lambda_m(C) = c_0 + \sum_{k=1}^{n/2-1} 2c_k\cos(2\pi mk/n) + c_{n/2}\cos (m\pi), \quad m = 0,\ldots,n-1.
\end{equation}

The circulant matrices form a commutative matrix $*$-algebra, as do the symmetric circulant matrices.
In particular, all circulant matrices share a set of eigenvectors, given by
the columns of the \emph{discrete Fourier transform matrix}:
\begin{equation}
\label{def:discrete Fourier transform matrix}
Q_{ij} := \frac{1}{\sqrt{n}}e^{-2\pi \sqrt{-1}ij/n}, \quad i,j = 0,\ldots,n-1.
\end{equation}
One has $Q^*Q = I$, and
 $Q^*CQ$ is a diagonal matrix for any circulant matrix $C$.
Also note that $Q^*e = \sqrt{n}e$.

We may define a basis $C^{(0)},\ldots,C^{\lfloor n/2\rfloor }$ for the symmetric circulant matrices as follows:
to obtain $C^{(i)}$ we set $c_i = c_{n-i} = 1$ in (\ref{def:circulant matrix}) and all other $c_j$'s to zero.

(We set $C_0 = 2I$ and also multiply $C_{n/2}$ by $2$ if $n$ is even).

By (\ref{eigs symmetric circulant}) and (\ref{eigs symmetric circulant even case}),
 the eigenvalues of these basis matrices are:
\begin{equation}
\lambda_m(C^{(k)}) =  2\cos(2\pi mk/n), \quad m = 0,\ldots,n-1, \; k = 0,\ldots,\lfloor n/2 \rfloor.
\label{eigs of circulant basis}
\end{equation}
Also note that
\[
\lambda_m(C^{(k)}) = \lambda_{n-m}(C^{(k)}), \quad m = 1,\ldots, \lfloor n/2 \rfloor, \; k = 0,\ldots,\lfloor n/2 \rfloor
\]
so that each matrix $C^{(k)}$ has only $1+ \lfloor n/2 \rfloor$ distinct eigenvalues.

\end{document}